\documentclass[MSNbibl,number,citesort,dvips]{arxbj}
\usepackage{mathrsfs}

\aid{0}
\volume{19}
\issue{1}
\pubyear{2013}
\firstpage{1}
\lastpage{17}
\doi{10.3150/11-BEJ389}

\makeatletter
\newcommand{\mv}{\mathrm{E}}
\newcommand{\diag}{\operatorname{diag}}
\newcommand{\df}{\mathrm{d}}
\newcommand{\e}{\mathrm{e}}
\newcommand{\Nset}{\mathbb{N}}
\newcommand{\Rset}{\mathbb{R}}
\newcommand{\Cset}{\mathbb{C}}
\newcommand{\Eop}{\mathscr{E}}
\newcommand{\Kop}{\mathscr{K}}
\newcommand{\Lop}{\mathscr{L}}
\newcommand{\Gop}{\mathscr{G}}
\newcommand{\Mop}{\mathscr{M}}
\newcommand{\Iop}{\mathbb{I}}
\newcommand{\Green}{\mathcal{G}}

\newtheorem{theo}{Theorem}
\newtheorem{prop}{Proposition}
\newtheorem{corr}{Corollary}
\newremark{remark}{Remark}
\newremark{example}{Example}

\renewcommand{\epsilon}{\varepsilon}
\newcommand{\eqref}[1]{(\ref{#1})}
\newcommand{\fracc}[2]{{#1}/{(#2)}}
\newcommand{\fracf}[2]{({#1})/({#2})}
\makeatother

\begin{document}
\begin{frontmatter}

\title{Functional data analysis in an operator-based mixed-model framework}
\runtitle{Mixed models for functional data}

\begin{aug}
\author{\fnms{Bo} \snm{Markussen}\corref{}\ead[label=e1]{bomar@life.ku.dk}}
\address{Department of Basic Sciences and Environment,
  University of Copenhagen, 1871 Frederiksberg~C, Denmark. \printead{e1}}
\end{aug}

\received{\smonth{11} \syear{2010}}
\revised{\smonth{7} \syear{2011}}

\begin{abstract}
Functional data analysis in a mixed-effects model framework is done
using operator calculus. In this approach the functional parameters are
treated as serially correlated effects giving an alternative to the
penalized likelihood approach, where the functional parameters are
treated as fixed effects. Operator approximations for the necessary
matrix computations are proposed, and semi-explicit and numerically
stable formulae of linear computational complexity are derived for
likelihood analysis. The operator approach renders the usage of a
functional basis unnecessary and clarifies the role of the boundary
conditions.
\end{abstract}


\end{frontmatter}

\section{Introduction} \label{secintro}

The aim of this paper is to derive operator approximations of the
matrix computations used to estimate the fixed and the random effects
in a mixed-effects model, where $M$ samples $y_1,\ldots,y_M \in
\Rset^N$ of temporal curves have been observed at $N$ predefined time
points $t_1,\ldots,t_N$. The main technical contribution of this
paper, making it practically possible to solve the estimation problem
as a functional estimation problem, is that the proposed operator
approximations have linear computational complexity in the sample
length $N$. Consequently, the mixed-effects inference becomes feasible
in the realm of functional data analysis, where $N$ can be
large.

Concatenating the samples $y_m = \{ y_{mn} \}_{n=1,\ldots,N} \in
\Rset^N$ into an observation vector $y = \{ y_m \}_{m=1,\ldots,M} \in
\Rset^{N_{\mathrm{total}}}$ with dimension $N_{\mathrm{total}} = N \ast M$
the statistical model we use is given by
\begin{equation} \label{eqmodel}
  y = \Gamma\beta + Z u + x + \epsilon.
\end{equation}
In this linear mixed-effects model the design matrices
$\Gamma \in \Rset^{N_{\mathrm{total}} \times p}$ and
$Z \in \Rset^{N_{\mathrm{total}} \times q}$ are known and assumed to have
full ranks $p$ and $q$, respectively, and the fixed effects $\beta \in
\Rset^p$ and the random effects $u \sim \mathcal{N}_q(0,\sigma^2 G)$
may be shared by the $M$ samples. The random component\vadjust{\goodbreak} $x = \{ x_m
\}_{m=1,\ldots,M} \sim \mathcal{N}_{N_{\mathrm{total}}}(0,\sigma^2 R)$
is partitioned in the same way as the observation vector $y$ and
consists of discretized readings $x_m = \{ x^{\mathrm{fct}}_m(t_n)
\}_{n=1,\ldots,N} \in \Rset^N$ of unobserved (latent) random functions
$x^{\mathrm{fct}}_m\dvtx [a,b] \to \Rset$. We assume that the random
functions $x^{\mathrm{fct}}_1,\ldots,x^{\mathrm{fct}}_M$ are independent
and identically distributed Gaussian processes with zero mean. The
covariance matrix $\sigma^2 R$ will be specified below appealing to
the smoothing splines methodology often used in functional data
analysis. Due to the i.i.d. assumption there exists a covariance
matrix $R_0 \in \Rset^{N \times N}$ such that $R = R_0 \otimes
\Iop_M$, where $\otimes$ is the Kronecker tensor product, and $\Iop_M
\in \Rset^{M \times M}$ is the identity matrix of dimension $M$. The
last component in the mixed-effects model is the measurement noise
$\epsilon \sim \mathcal{N}_{N_{\mathrm{total}}}(0,\sigma^2
\Iop_{N_{\mathrm{total}}})$.

Our objective is to derive computationally efficient formulae for the
maximum likelihood estimate of the fixed effects $\beta$, the best
linear unbiased predictions of the random effects $u \in \Rset^q$ and
of the latent random functions $x^{\mathrm{fct}}_m\dvtx [a,b] \to \Rset$,
and for the restricted likelihood function. The latter allow for
restricted likelihood inference on the variance parameters $\sigma^2 >
0$, $G \in \Rset^{q \times q}$ and $R_0 \in \Rset^{N \times N}$. The
methodology presented in this paper has two notable differences as
compared to the penalized likelihood approach to functional data
analysis; see, for example, the books by Ramsay and Silverman
\cite{RamsaySilverman02,RamsaySilverman05}. From the viewpoint of
computations we devise methods that work directly on the data
vector $y$ and, for example, provide predictions $\partial_t^\mu
\mv[x^{\mathrm{fct}}_m(t) \vert y]$ of the temporal derivatives of the
latent functional parameters. In particular, there is no basis
representation of the functional object $\mv[x^{\mathrm{fct}}_m \vert
y]$. This is by contrast with the standard technology used in functional
data analysis, where functional parameters are given a finite
dimensional representation, for example,  in a spline basis, and the
sparseness of the associated covariance matrices is invoked to achieve
feasible computations. As an alternative to this we use analytically
tractable operator approximations of the matrix equations. From the
viewpoint of statistical modeling we model the functional parameters
$x^{\mathrm{fct}}_m$ as random effects. Whether this is preferable over
the fixed effect interpretation underlying the penalized likelihood
depends on the particular application at hand. The distinction between
random and fixed effects is here the same as for classical
mixed-effects models; see \cite{Robinson91} for a thorough discussion
of the issue and \cite{KouEfron02} for a comparison of the associated
inference methodologies.

In the simplified version $y = x + \epsilon$ of model
equation \eqref{eqmodel}, the sample $y_m$ may be understood as a noisy
observation of the function $x^{\mathrm{fct}}_m\dvtx [a,b] \to \Rset$
taken at the sample points $t_1,\ldots,t_N$. In the penalized
likelihood approach to functional data analysis the functional
parameters $x^{\mathrm{fct}}_m$ are treated as fixed effects. The
penalized negative log likelihood is given by
\begin{equation} \label{eqnlogpenalizedlikelihood}
  N_{\mathrm{total}} \log\sigma + \frac{1}{2\sigma^2} \sum_{m=1}^M
   \Biggl( \sum_{n=1}^N | y_{mn} -  x^{\mathrm{fct}}_m(t_n) |^2
  + \lambda \int_a^b | \Kop x^{\mathrm{fct}}_m(t) |^2\, \df t
   \Biggr),
\end{equation}
where $\Kop$ is a differential operator of some order $k$ measuring the
roughness of a function $\theta \in \mathcal{C}^k([a,b] ; \Rset)$. The
so-called smoothing parameter $\lambda > 0$ quantifies the trade-off
between a close fit of the observations and the roughness of the
functional parameters. Since the space of functions is
infinite-dimensional, such a trade-off is required to avoid overfitting
of the finite number of data points.

In this paper we avoid the curse of dimensionality by providing the
theoretical solution in the function space before plugging in the
observed grid readings to compute the solution.\vadjust{\goodbreak} This is done using the
operator $\Lop = \Kop^\dag \Kop$, which is of order $2k$ and defined
on $\mathcal{C}^{2k}([a,b] ; \Rset)$. To ensure positive definiteness
of $\Lop$ we impose boundary conditions. Let $a_i, b_i \in
\{i-1,2k-i\}$ for $i=1,\ldots,k$ be fixed, and let the function space
$\mathcal{H}$ be defined by
\begin{equation} \label{defH}
  \mathcal{H} = \bigl \{ \theta \in \mathcal{C}^{2k} ([a,b] ; \Rset )
   \vert \theta^{(a_i)}(a) = \theta^{(b_i)}(b) = 0
  \mbox{ for $i=1,\ldots,k$}  \bigr\},
\end{equation}
where $\theta^{(i)}$ denotes the $i$th order derivative of
$\theta$. Applying integration by parts $k$ times the penalty terms in
equation \eqref{eqnlogpenalizedlikelihood} may be rewritten via
\[
  \int_a^b | \Kop \theta(t) |^2\, \df t
= \int_a^b \theta(t) \Lop \theta(t)\, \df t, \qquad \theta \in \mathcal{H}.
\]
This identity also implies that $\Lop$ is a positive semidefinite operator
on $\mathcal{H}$. A condition ensuring $\Lop$ to be invertible is
given in Section~\ref{secexplicit}. In the affirmative case the
inverse operator is given by a so-called Green's function $\Green(t,s)$
via $\Lop^{-1} f(t) = \int_a^b \Green(t,s) f(s)\, \df s$. Since $\Lop$
is positive definite it follows that $\Green(t,s)$ is positive
definite. In particular, the matrix defined by
\[
R_0 =  \{ \Green(t_n,t_m)  \}_{n,m=1,\ldots,N} \in \Rset^{N
  \times N}
\]
is positive definite and may be used as the variance of the serially
correlated effects $x_m$. This specification establishes a link
between the covariance matrix $\sigma^2 R_0$ of the discretized
readings $x_m$ in the model equation~\eqref{eqmodel} and the penalized
likelihood equation \eqref{eqnlogpenalizedlikelihood}.

The proposed methodology can be slightly generalized taking $\Lop$ as
the sum of squares $\sum_{l=1}^L \Kop_l^\dag \Kop_l$ of operators
measuring different aspects of roughness. The operator $\Lop$ may be
interpreted as a precision and used in the parameterization of a
statistical model. This is by contrast with standard software for
mixed-effects models such as the nlme-package \cite{PinheiroBates00}
in R or the MIXED procedure in SAS, where the parameterization is done
in terms of variances. In \cite{Verbylaetal99} a similar approach was
taken for the analysis of longitudinal data, and further references
may be found in \cite{LeeNelderPawitan06}, Chapter 8.4.

The remainder of this paper is organized as follows.
Section~\ref{secinference} reviews inference techniques for the
model equation \eqref{eqmodel}. In particular, we present the matrix
formulae that will be approximated by their operator equivalents.
Section~\ref{secmath} provides the mathematical contributions of
the paper. In this section the operator approximation is introduced and
refined for the case of equidistant observations, that is, $t_n = a +
\frac{2n-1}{N}(b-a)$. In particular, we derive semi-explicit and
numerically stable formulae for the needed computations in the case of
equidistant observations. In Section~\ref{secstatistics} the operator
approximation is applied on the matrix formulae from
Section~\ref{secinference}. This leads to concrete algorithms that
have been implemented in an R-package named fdaMixed \cite{Markussen11}.

\section{Inference in the mixed-effects model}
\label{secinference}

This section reviews estimation and inference techniques for the model
equation \eqref{eqmodel}. Since the derivations of the matrix formulae
stated below are standard (see,
e.g., \cite{Harville77,Henderson75,HendersonSearle81,Robinson91}),
no proofs will be given. The dimensions are given by
\begin{eqnarray*}
  y  = \Gamma\beta + Z u + x + \epsilon \in \Rset^{N_{\mathrm{total}}},
  \qquad
  \beta  \in \Rset^p,
  u  \in \Rset^q,
  x  \in \Rset^{N_{\mathrm{total}}},
  \epsilon &\in \Rset^{N_{\mathrm{total}}},
\end{eqnarray*}
where $N_{\mathrm{total}} = N \ast M$. Based on the covariance matrices
$G$ and $R = R_0 \otimes \Iop_M$ we define the matrices $A_0 = \Iop_N
+ R_0$, $A = \Iop_{N_{\mathrm{total}}} + R = A_0 \otimes \Iop_M$, and
\begin{eqnarray*}
  C_u  =  ( G^{-1} + Z^\top A^{-1} Z  )^{-1},  \qquad
  C_r  = A^{-1} - A^{-1} Z C_u Z^\top A^{-1},  \qquad
  C_\beta  =  ( \Gamma^\top C_r \Gamma  )^{-1}.
\end{eqnarray*}
The matrix formulae will be stated such that for moderately sized $p$ and
$q$ the computational obstacle of their practical implementation lies
in the initialization and inversion of the $N$-dimensional matrix
$A_0$. The circumvention of this obstacle is the topic of
Section~\ref{secmath}.

For known variance parameters $\sigma^2$, $G$, $R_0$, the best unbiased
estimate for the fixed effects is given by the maximum likelihood estimate
\begin{equation} \label{eqMLE}
  \hat{\beta} = C_\beta \Gamma^\top C_r y = C_\beta \Gamma^\top  (
  A^{-1} y - A^{-1} Z C_u Z^\top A^{-1} y  ).
\end{equation}
The best linear unbiased predictions (BLUPs) for the random effects $u$
and the serially correlated effects $x = \{ x_m \}_{m=1,\ldots,M}$ are
given by the conditional means
\begin{equation} \label{eqBLUP}
  \mv[u \vert y]  = C_u Z^\top A^{-1}  ( y - \Gamma\hat{\beta}  ),
  \qquad
  \mv[x \vert y]  = R A^{-1}  ( y - \Gamma\hat{\beta} -
  Z \mv[u \vert y]  ).
\end{equation}

It is generally agreed (see, e.g., \cite{Harville77} and
\cite{LeeNelderPawitan06}, Chapter 5.3) that the variance parameters
may be estimated as the maximizers of the restricted likelihood. One
of the factors in the likelihood is the determinant of $A_0 = \Iop_N +
R_0$. To derive the operator approximation of this factor we use the
representation
\begin{equation} \label{eqlogdet}
  \log\det A_0 = \int_0^1 \partial_v \log\det(\Iop + v R_0)\,
  \df v = \int_0^1 \sum_{j=1}^N e_j^\top  ( v \Iop_N + R_0^{-1}
   )^{-1} e_j\, \df v,
\end{equation}
where the vectors $e_j = \{ 1_{j=n} \}_{n=1,\ldots,N} \in \Rset^N$ for
$j=1,\ldots,N$ constitute an orthonormal basis for
$\Rset^N$. Using this representation and introducing the
conditional residuals $r = y - \Gamma \hat{\beta} - Z \mv[u \vert y] -
\mv[x \vert y]$, the double negative log restricted likelihood is
given by
\begin{eqnarray} \label{eqnlogRelikelihood}
  &&(2N_{\mathrm{total}}-2p) \log\sigma
+ M \int_0^1 \sum_{j=1}^N e_j^\top  ( v \Iop_N + R_0^{-1}  )^{-1}
e_j\, \df v \nonumber\\
&& \quad {}+ \log\det ( \Iop_q + Z^\top A^{-1} Z G  )
- \log\det C_\beta \\
&& \quad {}+ \sigma^{-2}  ( r^\top r + \mv[u \vert y]^\top G^{-1} \mv[u \vert y]
+ \mv[x \vert y]^\top R^{-1} \mv[x \vert y]  ),\nonumber
\end{eqnarray}
where it should be kept in mind that $C_\beta$, $r$, $\mv[u \vert y]$,
$\mv[x \vert y]$ depend on $G$ and $R_0$. The profile estimate
for the error variance $\sigma^2$ has an explicit form,
\[
  \hat{\sigma}^2 = \frac{1}{N_{\mathrm{total}}-p}  ( r^\top r
+ \mv[u \vert y]^\top G^{-1} \mv[u \vert y]
+ \mv[x \vert y]^\top R^{-1} \mv[x \vert y]  ).
\]

We conclude this section by reviewing some theoretical results on the
inference techniques described above. The errors $\hat{\beta}-\beta$,
$\mv[u \vert y]-u$, $\mv[x \vert y]-x$ follow a joint Gaussian
distribution, and their joint covariance may be derived using
\cite{Henderson75}, Section 2.4. Kackar and Harville
\cite{KackarHarville81} show that if the estimators for the variance
parameters are translation-invariant and even functions of $y$, then
$\hat{\beta}$, $\mv[u \vert y]$, $\mv[x \vert y]$ remain unbiased when
the estimates are inserted in place of the unknown variance
parameters. As explained by Welham and
Thompson \cite{WelhamThompson97} inference on $\beta$ may be done as
$\chi^2$-tests on twice the log ratio between the maximum restricted
likelihoods, where the design matrix under the null hypothesis is used
in the definition of the restricted likelihood under the
model. Simulation studies done by Morrell \cite{Morrell98} suggest
that inference on the variance parameters may be done as
$\chi^2$-tests on twice the log ratio between the maximum restricted
likelihoods, but here the formal asymptotic theory appears to be less
developed.

\section{Functional embedding of discrete data} \label{secmath}

Functional data consist of observations of continuous curves at
discrete sample points. As an alternative to computations based on
spline representations and sparse matrix computations we embed the
discrete observations into the continuous setting and approximate the
matrix computations by their operator counterparts. In order to
describe this operator approximation we first introduce some notation.

By a discretization of size $N$ of the time interval $[a,b]$ we mean a
set of points $\mathcal{T} = \{t_1,\ldots,t_N\}$ with $a < t_1 <
\cdots < t_N < b$. Such a discretization is said to be equidistant if
$t_n = a + \frac{2n-1}{2N} (b-a)$, and in that case we associate the
mesh length given by $\Delta = (b-a)/N$. To ease notation we
implicitly adjoin the points $t_0=a$ and $t_{N+1}=b$ to any
discretization of size $N$.

Given a vector $z = \{z_n\}_{n=1,\ldots,N} \in \Rset^N$ we denote by
$\Eop_z$ the piecewise linear embedding of $z \in \Rset^N$ into
$\mathcal{C}([a,b] ; \Rset)$, that is, the function that is linear on the
segments $[t_n,t_{n+1}]$ for $n=0,\ldots,N$ with $\Eop_z(a)=z_1$,
$\Eop_z(b)=z_N$ and $\Eop_z(t_n) = z_n$ for $n=1,\ldots,N$. We also
introduce the multiplication operator $\Mop_{\mathcal{T}}$ on
$\mathcal{C}([a,b] ; \Rset)$ defined by
\[
  \Mop_{\mathcal{T}} f(t) = \Eop_\mu(t) f(t), \qquad f \in
  \mathcal{C}([a,b] ; \Rset),
\]
where $\mu = \{\mu_n\}_{n=1,\ldots,N} \in \Rset^N$ is given from the
discretization $\mathcal{T}$ via
\begin{equation} \label{eqmu}
\mu_n = \cases{\displaystyle
  2/(t_2+t_1-2a), &  \quad \mbox{for $n=1$,} \cr\displaystyle
  2/(t_{n+1}-t_{n-1}), & \quad  \mbox{for $n=2,\ldots,N-1$,} \cr\displaystyle
  2/(2b-t_N-t_{N-1}), & \quad  \mbox{for $n=N$.}
}
\end{equation}
In particular, if $\mathcal{T}$ is equidistant, then
$\Mop_{\mathcal{T}} = \Delta^{-1} \Iop$.

\begin{prop} \label{propapproximation}
Let a discretization $\mathcal{T}$ of the interval $[a,b]$, $t \in
[a,b]$ and $\Green \in \mathcal{C}([a,b] \times [a,b] ; \Rset)$ be given.
Assume that $\Green(t,\cdot)$ is twice differentiable on the segments
$[t_n,t_{n+1}]$ with continuous derivatives $\Green^{(i)}(t,\cdot)$. For
$z \in \Rset^N$ there exists $\xi_n \in (t_n,t_{n+1})$ for
$n=0,\ldots,N$ and $\zeta_1 \in (a,t_1)$, $\zeta_N \in (t_N,b)$ such
that $\sum_{n=1}^N \Green(t,t_n) z_n - \int_a^b \Green(t,s) \Eop_\mu(s)
\Eop_z(s)\, \df s$ equals
\begin{eqnarray} \label{eqerror}
&&\frac{(t_1-a)^2}{2} \Green^{(1)}(t,\zeta_1) \mu_1 z_1
- \frac{(b-t_N)^2}{2} \Green^{(1)}(t,\zeta_N) \mu_N z_N \nonumber\\
&& \quad {}+ \frac{1}{12} \sum_{n=0}^N (t_{n+1}-t_n)^3
 \biggl( \Green^{(2)}(t,\xi_n) \Eop_\mu(\xi_n) \Eop_z(\xi_n) +
  2 \Green^{(1)}(t,\xi_n) \frac{\mu_{n+1}-\mu_n}{t_{n+1}-t_n}
  \Eop_z(\xi_n) \\
&&  \hphantom{\quad {}+\frac{1}{12} \sum_{n=0}^N (t_{n+1}-t_n)^3
 \biggl(} {}+ 2 \Green^{(1)}(t,\xi_n) \Eop_\mu(\xi_n)\frac{z_{n+1}-z_n}{t_{n+1}-t_n}
+ 2 \Green(t,\xi_n) \frac{\mu_{n+1}-\mu_n}{t_{n+1}-t_n}
\frac{z_{n+1}-z_n}{t_{n+1}-t_n}  \biggr),\nonumber
\end{eqnarray}
where $\mu$ is given by equation~\eqref{eqmu}.
\end{prop}

\begin{pf}
The \textit{trapezoidal rule} of integration
\cite{KincaidCheney91}, Section 7.2, gives intermediate points $\xi_n
\in (t_n,t_{n+1})$ such that $\int_a^b \Green(t,s) \Eop_\mu(s)
\Eop_z(s)\, \df s$ equals\vspace*{-2pt}
\begin{eqnarray*}
&&\frac{t_1-a}{2} \Green(t,a) \Eop_\mu(a) \Eop_z(a)
+ \sum_{n=1}^N \frac{t_{n+1}-t_{n-1}}{2} \Green(t,t_n) \Eop_\mu(t_n)
\Eop_z(t_n) \\[-2pt]
&& \quad {}+ \frac{b-t_N}{2} \Green(t,b) \Eop_\mu(b) \Eop_z(b)
- \frac{1}{12} \sum_{n=0}^N (t_{n+1}-t_n)^3  ( \Green(t,\cdot)
  \Eop_\mu \Eop_z  )^{(2)}(\xi_n).
\end{eqnarray*}
The result follows inserting the piecewise linear functions $\Eop_\mu$
and $\Eop_z$, the first-order Taylor expansions at some
intermidiate points $\zeta_1 \in (a,t_1)$, $\zeta_N \in (t_N,b)$,
\begin{eqnarray*}
  \frac{t_1-a}{2} \Green(t,a) \Eop_\mu(a) \Eop_z(a)
&=& \frac{t_1-a}{2} \Green(t,t_1) \mu_1 z_1
- \frac{(t_1-a)^2}{2} \Green^{(1)}(t,\zeta_1) \mu_1 z_1, \\
  \frac{b-t_N}{2} \Green(t,b) \Eop_\mu(b) \Eop_z(b)
&=& \frac{b-t_N}{2} \Green(t,t_N) \mu_N z_N
+ \frac{(b-t_N)^2}{2} \Green^{(1)}(t,\zeta_N) \mu_N z_N,
\end{eqnarray*}
by expanding the second-order derivative and by rearranging the terms.
\end{pf}

\begin{corr} \label{corrapproximation}
If the discretization $\mathcal{T}$ is equidistant, then there exists
$\widetilde{\xi}_n \in (\xi_{n-1},\xi_n) \subset (t_{n-1},t_{n+1})$
for $n=1,\ldots,N$ such that the approximation error
equation \eqref{eqerror} equals\vspace*{-2pt}
\begin{eqnarray} \label{eqerror2}
&&\frac{b-a}{12 N} \sum_{n=1}^N  \bigl( t_{n+1}-t_{n-1} - 3
(\xi_n-\xi_{n-1})  \bigr) \Green^{(2)}(t,t_n) z_n \nonumber\\[-2pt]
&& \quad {}+ \frac{b-a}{12 N} \sum_{n=1}^N (t_{n+1}-\xi_n)
 \bigl( \Green^{(2)}(t,\xi_n) - \Green^{(2)}(t,t_n)  \bigr) z_n
 \nonumber\\[-2pt]
&& \quad {}- \frac{b-a}{6 N} \sum_{n=1}^N (\xi_n-\xi_{n-1})
 \bigl( \Green^{(2)}(t,\widetilde{\xi}_n) - \Green^{(2)}(t,t_n)
 \bigr) z_n \\[-2pt]
&& \quad {}+ \frac{b-a}{12 N} \sum_{n=1}^N (\xi_{n-1}-t_{n-1})
 \bigl( \Green^{(2)}(t,\xi_{n-1}) - \Green^{(2)}(t,t_n)  \bigr) z_n \nonumber\\[-2pt]
&& \quad {}+ \frac{(b-a)^2}{8 N^2} \Green^{(1)}(t,\zeta_1) \mu_1 z_1
- \frac{(b-a)^2}{8 N^2} \Green^{(1)}(t,\zeta_N) \mu_N z_N.
\nonumber
\end{eqnarray}
\end{corr}

\begin{pf}
Equidistant spacing implies that the factors $\mu_n = N/(b-a)$ defined
in equation \eqref{eqmu} are constant, and the approximation error
equation \eqref{eqerror} reduces to
\begin{eqnarray*}
&&\frac{(b-a)^2}{8 N^2} \Green^{(1)}(t,\zeta_1) \mu_1 z_1
- \frac{(b-a)^2}{8 N^2} \Green^{(1)}(t,\zeta_N) \mu_N z_N \\[-2pt]
&& \quad {}+ \sum_{n=0}^N \frac{b-a}{12 N} \biggl ( \frac{b-a}{N}
  \Green^{(2)}(t,\xi_n) \Eop_z(\xi_n) +
  2 \Green^{(1)}(t,\xi_n)  (z_{n+1}-z_n  )  \biggr).
\end{eqnarray*}
The last sum equals\vspace*{-1.5pt}
\begin{eqnarray*}
&&\sum_{n=0}^N \frac{b-a}{12 N} \Green^{(2)}(t,\xi_n)
 \bigl( (\xi_n-t_n) z_{n+1}  + (t_{n+1}-\xi_n) z_n  \bigr) \\
&& \quad {}+ \sum_{n=1}^N \frac{b-a}{6 N} \bigl (
\Green^{(1)}(t,\xi_{n-1}) - \Green^{(1)}(t,\xi_n)  \bigr) z_n.
\end{eqnarray*}
By Taylor's theorem there exists $\widetilde{\xi}_n \in
(\xi_{n-1},\xi_n)$ such that this equals
\begin{eqnarray*}
&&\sum_{n=1}^N \frac{b-a}{12 N}  \bigl( (t_{n+1}-\xi_n) \Green^{(2)}(t,\xi_n)
+ (\xi_{n-1}-t_{n-1}) \Green^{(2)}(t,\xi_{n-1}) \\
&& \hphantom{\sum_{n=1}^N \frac{b-a}{12 N}  \bigl(} {}- 2 (\xi_n-\xi_{n-1}) \Green^{(2)}(t,\widetilde{\xi}_n)  \bigr) z_n.
\end{eqnarray*}
The corollary follows centering the terms $\Green^{(2)}(t,\cdot)$
around $\Green^{(2)}(t,t_n)$.
\end{pf}

If the matrix $D \in \Rset^{N \times N}$ and the integral operator
$\Gop$ on $\mathcal{C}([a,b] ; \Rset)$ are defined by $D = \{
\Green(t_n,t_m) \}_{n,m=1,\ldots,N}$ and $\Gop f(t) = \int_a^b
\Green(t,s) f(s)\, \df s$, then the preceding results suggest the
approximation\vspace*{-1.5pt}
\begin{equation} \label{eqapprox}
  D z \approx  \{ \Gop \Mop_{\mathcal{T}} \Eop_z(t_n)
   \}_{n=1,\ldots,N} \in \Rset^N, \qquad z \in \Rset^N.
\end{equation}
Green's functions usually possess sufficient smoothness for
Proposition~\ref{propapproximation} to apply (see,
e.g., \cite{VillaneuvaJodar91}), and hence the approximation error in
equation \eqref{eqapprox} vanishes as $\max_{n=0,\ldots,N} | t_{n+1}
- t_n |$ goes to zero. In case of equidistant discretizations
this property is refined in Corollary~\ref{corrapproximation}. The
first term in equation \eqref{eqerror2} is of size $\mathrm{O}(N^{-1})$ and the
other terms are of size $\mathrm{O}(N^{-2})$. Perhaps the first term can be
used to derive and correct a bias arising from the proposed operator
approximation, but we will leave this to be studied in future work.

\subsection{Explicit operator computations} \label{secexplicit}

To motivate the derivations done in this section we may consider the
model equation \eqref{eqmodel} without the fixed and the random effects,
that is, $y=x+\epsilon$. In this case equation \eqref{eqapprox} implies the\vadjust{\goodbreak}
approximation of the prediction equation \eqref{eqBLUP} of the $m$th
serially correlated effect given by
\[
  \mv[ x_m \vert y] = R_0 A_0^{-1} y_m =  ( \Iop_N + R_0^{-1}
   )^{-1} y_m \approx  \{  ( \Iop + \Mop_{\mathcal{T}}^{-1}
  \Lop  )^{-1} \Eop_{y_m}(t_n)  \}_{n=1,\ldots,N}
\]
and the approximation of the logarithmic determinant
equation \eqref{eqlogdet} given by
\begin{eqnarray} \label{eqapprox2}
\log\det A_0 &=& \int_0^1 \sum_{j=1}^N e_j^\top  ( v \Iop_N + R_0^{-1}
 )^{-1} e_j\, \df v \nonumber
 \\[-9pt]
 \\[-9pt]
&\approx& \int_0^1 \sum_{j=1}^N  ( v\Iop + \Mop_{\mathcal{T}}^{-1} \Lop
 )^{-1} \Eop_{e_j}(t_j)\, \df v
\approx \int_0^1 \int_a^b \Green_v(t,t)\, \df t\, \df v,
\nonumber
\end{eqnarray}
where $\Green_v$ is the Green's function for $v\Iop +
\Mop_{\mathcal{T}}^{-1} \Lop$. As shown in
Section~\ref{secstatistics} the matrix formulae used for inference in
the full mixed-effects model equation \eqref{eqmodel} may be similarly
approximated. In order to develop our computational methodology we
derive semi-explicit and numerically stable inversion formulae for
differential operators of the type $\Lop_* = \Iop +
\Mop_{\mathcal{T}}^{-1} \Lop$. If the discretization $\mathcal{T}$ is
equidistant with mesh length $\Delta$, and the differential operator
$\Lop$ has constant coefficients, then $\Lop_* = \Iop + \Delta \Lop$
may be inverted using Theorem~\ref{thmGreen} stated below. Boundary
conditions play an essential role in this theorem, and the reader may
want to refresh the definition of the space $\mathcal{H}$ given in
equation \eqref{defH}.\vspace*{-2pt}

\begin{theo} \label{thmGreen}
Consider a differential operator $\Lop_*$ on $\mathcal{H}$ given by
\begin{equation} \label{eqLop}
  \Lop_* \theta(t) = \alpha_{2k} \theta^{(2k)}(t) + \alpha_{2k-1}
  \theta^{(2k-1)}(t) + \cdots + \alpha_1 \theta^{(1)}(t) + \alpha_0
  \theta(t)
\end{equation}
with $\alpha_{2k} \neq 0$. Let $J = \diag(J_1,\ldots,J_p) \in
\Cset^{2k \times 2k}$, with $J_j \in \Cset^{k_j \times k_j}$, be the
Jordan canonical form of the companion matrix
\begin{equation} \label{eqcompanion}
  C = \pmatrix{
0 & 1 & 0 & \cdots & 0 \cr
\vdots &\ddots & \ddots & \ddots & \vdots \cr
0 & \cdots & 0 & 1 & 0 \vspace*{4pt}\cr
-\displaystyle\frac{\alpha_0}{\alpha_{2k}} & -\displaystyle\frac{\alpha_1}{\alpha_{2k}} & \cdots &
\cdots &\displaystyle -\frac{\alpha_{2k-1}}{\alpha_{2k}}
} \in \Rset^{2k \times 2k}.
\end{equation}
Let $M \in \Cset^{2k \times 2k}$ be a non-trivial solution of the
matrix equation $C M = M J$, and let $M_{1j} \in \Cset^{1 \times k_j}$
be the decomposition of the first row of $M$ along the Jordan blocks
$J_j$. Let $v_1 = ( 1 \enskip \cdots \enskip 1 ) \in
\Rset^{1 \times k}$, $v_2 = ( 0 \enskip \cdots \enskip 0 \enskip 1
)^\top \in \Rset^{2k \times 1}$, and let $F_a, F_b \in
\Rset^{k \times 2k}$ be given by
\begin{eqnarray*}
  F_a  =   \{ 1_{j-1=a_i}  \}_{
  {\fontsize{7.6pt}{9pt}\selectfont{\begin{array}{l} i=1,\ldots,k\\  j=1,\ldots,2k \end{array}}}},  \qquad
  F_b  =   \{ 1_{j-1=b_i}  \}_{
  {\fontsize{7.6pt}{9pt}\selectfont{\begin{array}{l} i=1,\ldots,k\\ j=1,\ldots,2k \end{array}}}}.
\end{eqnarray*}
Let $\bar{v}_1 = ( v_1 \enskip v_1 ) \in \Rset^{1
  \times 2k}$, and let
$\bar{F}_a, \bar{F}_b, W \in \Rset^{2k \times 2k}$ be defined by
\begin{eqnarray*}
  \bar{F}_a = \pmatrix{ F_a \cr 0_{k \times 2k}},  \qquad
  \bar{F}_b = \pmatrix{ 0_{k \times 2k} \cr F_b},  \qquad
  W = \pmatrix{ M_{11} & \cdots & M_{1p} \cr M_{11} J_1 & \cdots &
  M_{1p} J_p \cr \vdots & & \vdots \cr M_{11} J_1^{2k-1} & \cdots &
  M_{1p} J_p^{2k-1} }.\vadjust{\goodbreak}
\end{eqnarray*}
If the matrix $H = \bar{F}_a W \exp(a J) + \bar{F}_b W \exp(b J)$ is
invertible, then $\Lop_*$ is invertible. In the affirmative case the
inverse operator is an integral operator $\Lop_*^{-1} f(t) = \int_a^b
\Green_*(t,s) f(s)\, \df s$, where the Green's function $\Green_*$ is given
by
\begin{equation} \label{eqGreen}
  \Green_*(t,s) = \cases{\displaystyle
\alpha_{2k}^{-1} \bar{v}_1 \exp(t J) H^{-1} \bar{F}_a W \exp \bigl((a-s)
J \bigr) W^{-1} v_2, &  \quad \mbox{for $s \le t$,} \cr\displaystyle
-\alpha_{2k}^{-1} \bar{v}_1 \exp(t J) H^{-1} \bar{F}_b W
\exp \bigl((b-s) J \bigr) W^{-1} v_2, &  \quad \mbox{for $t \le s$.}
}
\end{equation}
\end{theo}

\begin{pf}
The proof follows specializing and condensing
\cite{VillaneuvaJodar91}, Theorem 3. The signs of~\cite{VillaneuvaJodar91},
equation (3.15), equation (3.24), should be changed due to a mistake of
sign in \cite{VillaneuvaJodar91}, equation (3.9). We allow for leading
coefficient $\alpha_{2k} \neq 1$ and have interchanged the indices $k$
and $p$ to align with the notation used in the present paper.
\end{pf}

Formula \eqref{eqGreen} is explicit and most satisfactory from a
theoretical point of view. But from a practical point of view the formula
can be numerically unstable since the exponentials $\exp(tJ)$,
$\exp((a-s)J)$ and $\exp((b-s)J)$ are weighted against similar
exponentials in the definition of the matrix $H$. Imposing symmetry of
the Jordan matrix it is, however, possible to remove the potential
numerical instabilities.

\begin{prop} \label{propnumericalstable}
Suppose that the characteristic polynomial
\begin{equation} \label{eqcharacteristic}
  \alpha_{2k} z^{2k} + \alpha_{2k-1} z^{2k-1} + \cdots + \alpha_1 z +
  \alpha_0 = 0
\end{equation}
for the differential operator \eqref{eqLop} has $2k$ distinct roots
$\eta^-_1,\eta^+_1,\ldots,\eta^-_k,\eta^+_k \in \Cset$
such that the real values of the $k$ eigenvalues
$\eta^-_1,\ldots,\eta^-_k$ are non-positive and the real values
of the $k$ eigenvalues $\eta^+_1,\ldots,\eta^+_k$ are
non-negative. Then the Jordan canonical form of the companion matrix
equation \eqref{eqcompanion} is diagonal with block diagonals consisting
of eigenvalues with non-positive and non-negative real values,
respectively,
\begin{eqnarray*}
J = \pmatrix{ J_- & 0_{k \times k} \cr 0_{k \times k} & J_+
},  \qquad
J_- = \diag(\eta^-_1,\ldots,\eta^-_k),  \qquad
J_+ = \diag(\eta^+_1,\ldots,\eta^+_k),
\end{eqnarray*}
and the matrix $W = { W_- \choose W_+ } \in
\Cset^{2k \times 2k}$ may be decomposed via
$W_-, W_+ \in \Cset^{2k \times k}$ defined by
\begin{eqnarray*}
  W_- = \pmatrix{
1 & \cdots & 1 \cr
\eta^-_1 & \cdots & \eta^-_k \cr
\vdots & \ddots & \vdots \cr
(\eta^-_1)^{2k-1} & \cdots & (\eta^-_k)^{2k-1}
},  \qquad
  W_+ = \pmatrix{
1 & \cdots & 1 \cr
\eta^+_1 & \cdots & \eta^+_k \cr
\vdots & \ddots & \vdots \cr
(\eta^+_1)^{2k-1} & \cdots & (\eta^+_k)^{2k-1}
}.
\end{eqnarray*}
Furthermore, define $v_1 = ( 1 \enskip \cdots \enskip 1 )
\in \Rset^{1 \times k}$, $v_2 = ( 0 \enskip \cdots \enskip 0 \enskip 1
  )^\top \in \Rset^{2k \times 1}$, $v_-, v_+ \in
  \Rset^{k \times 1}$ via
$W^{-1} v_2 = { v_- \choose v_+ }$, and the
vectors $\phi_\mu(t), \psi_\mu(t) \in \Rset^{1 \times k}$ for $t \in
[a,b]$ and $\mu \in \Nset_0$ by
\begin{eqnarray*}
  \phi_\mu(t) &=& \alpha_{2k}^{-1}  \bigl( v_1 J_-^\mu - v_1 J_+^\mu
  \e^{-(b-t) J_+} (F_b W_+)^{-1} F_b W_- \e^{(b-t) J_-}  \bigr) \\
&&{}\cdot \bigl ( I_{k \times k} - \e^{(t-a) J_-} (F_a W_-)^{-1} F_a W_+
\e^{-(b-a) J_+} (F_b W_+)^{-1} F_b W_- \e^{(b-t) J_-}  \bigr)^{-1}
\end{eqnarray*}
and
\begin{eqnarray*}
  \psi_\mu(t) &= &\alpha_{2k}^{-1} \bigl ( v_1 J_+^\mu - v_1 J_-^\mu
  \e^{(t-a) J_-} (F_a W_-)^{-1} F_a W_+ \e^{-(t-a) J_+}  \bigr) \\
&&{}\cdot  \bigl( I_{k \times k} - \e^{-(b-t) J_+} (F_b W_+)^{-1} F_b W_-
\e^{(b-a) J_-} (F_a W_-)^{-1} F_a W_+ \e^{-(t-a) J_+}  \bigr)^{-1}.
\end{eqnarray*}
Then the $\mu$th partial derivative $\partial_t^\mu \Green_*(t,s)$ of
the Green's function defined in equation \eqref{eqGreen} may be rewritten
as the numerically stable expression
\begin{equation} \label{eqGreennum1}
\cases{\displaystyle
  \phi_\mu(t) \e^{(t-s) J_-}  \bigl( v_- + \e^{(s-a) J_-} (F_a W_-)^{-1}
  F_a W_+ \e^{-(s-a) J_+} v_+  \bigr), & \quad  \mbox{for $s \le t$,} \cr\displaystyle
-\psi_\mu(t) \e^{-(s-t) J_+}  \bigl( v_+ + \e^{-(b-s) J_+} (F_b W_+)^{-1}
  F_b W_- \e^{(b-s) J_-} v_-  \bigr), & \quad  \mbox{for $t \le s$.}
}
\end{equation}
\end{prop}

\begin{pf}
From equation \eqref{eqGreen} we have that $\partial_t^\mu \Green_*(t,s)$
equals
\[
\cases{\displaystyle
\alpha_{2k}^{-1} \bar{v}_1 J^\mu \exp(t J) H^{-1} \bar{F}_a W
\exp \bigl((a-s) J \bigr) W^{-1} v_2 ,&  \quad \mbox{for $s \le t$,} \cr\displaystyle
-\alpha_{2k}^{-1} \bar{v}_1 J^\mu \exp(t J) H^{-1} \bar{F}_b W
\exp \bigl((b-s) J \bigr) W^{-1} v_2 ,&  \quad \mbox{for $t \le s$.}
}
\]
The crux of the reformulation of this representation lies in the
inversion of the matrix $H = \bar{F}_a W \exp(a J) + \bar{F}_b W
\exp(b J)$. To this end, we write $H \e^{-t J}$ and $\e^{t J} H^{-1}$
as block matrices with $k \times k$-blocks,
\begin{eqnarray*}
  H \e^{-t J} =\pmatrix{
F_a W_- \e^{(a-t) J_-} & F_a W_+ \e^{(a-t) J_+} \cr
F_b W_- \e^{(b-t) J_-} & F_b W_+ \e^{(b-t) J_+}
}, \quad
\e^{t J} H^{-1} = \pmatrix{
A_{11} & A_{12} \cr
A_{21} & A_{22}
}.
\end{eqnarray*}
Using elementary matrix algebra we find that $w_1 A_{11} + w_2 A_{21}$
for general $w_1, w_2 \in \Rset^{1 \times k}$ equals
\begin{eqnarray} \label{eqGreennum4}
 &\displaystyle \bigl( w_1 - w_2 \e^{-(b-t) J_+} (F_b W_+)^{-1} F_b W_- \e^{(b-t) J_-}
 \bigr)& \nonumber
 \\[-8pt]
 \\[-8pt]
 &\displaystyle \bigl( F_a W_- \e^{(a-t) J_-} - F_a W_+ \e^{-(b-a) J_+} (F_b W_+)^{-1}
  F_b W_- \e^{(b-t) J_-}  \bigr)^{-1}.&
\nonumber
\end{eqnarray}
Inserting this above we have that $\partial_t^\mu \Green_*(t,s)$ for
$s \le t$ equals
\[
\alpha_{2k}^{-1}
\pmatrix{ v_1 J_-^\mu & v_1 J_+^\mu }
\pmatrix{ A_{11} & A_{12} \cr A_{21} & A_{22} }
\pmatrix{ F_a W_- \e^{(a-s) J_-} & F_a W_+ \e^{(a-s) J_+} \cr
0_{k \times k} & 0_{k \times k} }
\pmatrix{ v_- \cr v_+ }
\]
which equals
\begin{equation} \label{eqGreennum5}
  \alpha_{2k}^{-1}  ( v_1 J_-^\mu A_{11} + v_1 J_+^\mu A_{21}  )
 \bigl( F_a W_- \e^{(a-s) J_-} v_- + F_a W_+ \e^{(a-s) J_+} v_+  \bigr).
\end{equation}
Combining equations \eqref{eqGreennum4} and \eqref{eqGreennum5} and
rearranging the exponential factors we arrive at
equation \eqref{eqGreennum1} for $s \le t$. The reformulation is done
similarly for $t \le s$.
\end{pf}

\begin{remark*}
From the viewpoint of statistical modeling, the results in
\cite{VillaneuvaJodar91} are more general in two valuable
ways. Firstly, the boundary conditions separately given at the
end-points of the sample interval via the matrices $F_a$, $F_b$ in
Theorem~\ref{thmGreen} may be given in form of linear combinations of
the curve and its derivatives at $a$ and $b$ via general $\bar{F}_a$
and $\bar{F}_b$. In particular, boundary conditions enforcing
periodicity may be stated. But to derive the numerically stable
formulae stated in Proposition~\ref{propnumericalstable}, we have
refrained from this possibility. Secondly, the results in
\cite{VillaneuvaJodar91} are given for matrix-valued functions. This
generalization allows our methods to be extended to multivariate
functional data analysis.
\end{remark*}

In the following theorem the explicit inversion formula is applied to
derive a simultaneous computation of $\partial_t^\mu (\Iop + \Delta
\Lop)^{-1} \Eop_z(t_n)$ for $n=1,\ldots,N$, where $z \in \Rset^N$,
that easily may be implemented with computational complexity
$\mathrm{O}(N)$. Furthermore, the inner integral in the approximation
equation \eqref{eqapprox2} of the logarithmic determinant may be
explicitly computed for Lebesgue almost all $v \in [0,1]$. In the
statement of the theorem we denote by $\odot$ the element-wise
multiplication of matrices or vectors of the same dimension. Unless
specified otherwise the $\odot$ operation is performed after ordinary
matrix multiplications.

\begin{theo} \label{thmfda}
Suppose the discretization $\mathcal{T}$ is equidistant with mesh
length $\Delta = (b-a)/N$, and assume that the operator in
equation \eqref{eqLop} given by $\Lop_* = \Iop + \Delta \Lop$ satisfies
the conditions of Proposition~\ref{propnumericalstable}. Denote by
$\Green_*$ the Green's function for $\Lop_*$, let $J_-$, $J_+$, $W_-$,
$W_+$, $v_-$, $v_+$, $\phi_\mu(t)$, $\psi_\mu(t)$ be as defined in
Proposition~\ref{propnumericalstable}, and let $\xi_-, \xi_-^0,
\xi_-^1, \xi_+, \xi_+^0, \xi_+^1 \in \Rset^{k \times 1}$ be defined by
\begin{eqnarray*}
  \xi_- &=& \biggl\{ \frac{\exp(\Delta \eta_i^- / 2) -
  1}{\eta_i^-} \biggr\}_{i=1,\ldots,k},  \hspace*{35.5pt}\qquad
  \xi_+ =\biggl\{ \frac{1 - \exp(-\Delta \eta_i^+ /
  2)}{\eta_i^+} \biggr\}_{i=1,\ldots,k}, \\[-2pt]
  \xi_-^0 &=& \biggl\{ \frac{1 - (1-\Delta\eta_i^-) \exp(\Delta
  \eta_i^-)}{\Delta (\eta_i^-)^2} \biggr\}_{i=1,\ldots,k}, \qquad
  \xi_+^0 = \biggl\{ \frac{\exp(-\Delta \eta_i^+) - 1 +
  \Delta \eta_i^+}{\Delta (\eta_i^+)^2}
  \biggr\}_{i=1,\ldots,k}, \\[-2pt]
  \xi_-^1 &=& \biggl\{ \frac{\exp(\Delta\eta_i^-) - 1 -
  \Delta\eta_i^-)}{\Delta (\eta_i^-)^2} \biggr\}_{i=1,\ldots,k},
  \hspace*{10pt}\qquad
  \xi_+^1 = \biggl\{ \frac{1 - (1+\Delta \eta_i^+)
  \exp(-\Delta \eta_i^+)}{\Delta (\eta_i^+)^2}
  \biggr\}_{i=1,\ldots,k}.
\end{eqnarray*}
For $z = \{ z_j \}_{j=1,\ldots,N} \in \Rset^N$ the $\mu$th derivative
$\partial_t^\mu (\Iop + \Delta \Lop)^{-1} \Eop_z(t_n)$ taken at the
sample point $t_n$ is given by
\begin{eqnarray*}
&&\phi_\mu(t_n) 1_{n > 1}  \sum_{j=1}^{n-1} \e^{(t_n-t_{j+1}) J_-} (v_-
\odot \xi_-^0) z_j \\[-2pt]
&&\quad{}+ \phi_\mu(t_n)  \sum_{j=1}^n \e^{(t_n-t_j) J_-}  \bigl(v_- \odot (1_{j=1}
\xi_- + 1_{j > 1} \xi_-^1)  \bigr) z_j  \\[-2pt]
&&\quad{}- \psi_\mu(t_n)  \sum_{j=n}^N \e^{-(t_j-t_n) J_+} \bigl (v_+ \odot
  (1_{j < N} \xi_+^0 + 1_{j=N} \xi_+)  \bigr) z_j \\[-2pt]
&&\quad{}- \psi_\mu(t_n) 1_{n < N}  \sum_{j=n+1}^N \e^{-(t_{j-1}-t_n) J_+} (v_+
  \odot \xi_+^1) z_j \\[-2pt]
&&\quad{}+ \phi_\mu(t_n) \e^{(t_n-a) J_-} (F_a W_-)^{-1} F_a W_+ 1_{n > 1}
 \sum_{j=1}^{n-1} \e^{-(t_j-a) J_+} (v_+ \odot \xi_+^0) z_j \\[-2pt]
&&\quad{}+ \phi_\mu(t_n) \e^{(t_n-a) J_-} (F_a W_-)^{-1} F_a W_+  \sum_{j=1}^n
\e^{-(t_{j-1}-a) J_+} \bigl (v_+ \odot (1_{j=1} \xi_+ + 1_{j > 1}
\xi_+^1)  \bigr) z_j \\[-2pt]
&&\quad{}- \psi_\mu(t_n) \e^{-(b-t_n) J_+} (F_b W_+)^{-1} F_b W_-  \sum_{j=n}^N
\e^{(b-t_{j+1}) J_-} \bigl (v_- \odot (1_{j < N} \xi_-^0 + 1_{j=N}
\xi_-)  \bigr) z_j \\[-2pt]
&&\quad{}- \psi_\mu(t_n) \e^{-(b-t_n) J_+} (F_b W_+)^{-1} F_b W_- 1_{n < N}
 \sum_{j=n+1}^N \e^{(b-t_j) J_-} (v_- \odot \xi_-^1) z_j.
\end{eqnarray*}

Concerning the log determinant assume that the operator in
equation \eqref{eqLop} given by $\Lop_* = v\Iop + \Delta \Lop$ for
fixed $v \in [0,1]$ satisfies the conditions of
Proposition~\ref{propnumericalstable}. Let the matrices $A_{--},
A_{++}, A_{-+}, A_{+-} \in \Rset^{k \times k}$ be defined by
\begin{eqnarray*}
  A_{--} &=& \biggl\{ 1_{i=j} N \e^{(b-a) \eta_i^-} +
  1_{i \neq j} \frac{\e^{(b-a) \eta_i^-} - \e^{(b-a)
  \eta_j^-}}{\Delta (\eta_i^- - \eta_j^-)}
  \biggr\}_{i,j=1,\ldots,k}, \\[-2pt]
  A_{++} &=& \biggl\{ 1_{i=j} N \e^{-(b-a) \eta_i^+}
  + 1_{i \neq j} \frac{\e^{-(b-a) \eta_i^+} - \e^{-(b-a)
  \eta_j^+}}{\Delta (-\eta_i^+ + \eta_j^+)}
  \biggr\}_{i,j=1,\ldots,k}, \\[-2pt]
  A_{-+} &=& \biggl\{ \frac{1 - \e^{(b-a)(\eta_i^- -
  \eta_j^+)}}{\Delta (-\eta_i^- + \eta_j^+)}
  \biggr\}_{i,j=1,\ldots,k}, \\[-2pt]
  A_{+-} &=& \biggl\{ \frac{1 - \e^{-(b-a)(\eta_i^+ -
  \eta_j^-)}}{\Delta (\eta_i^+ - \eta_j^-)}
  \biggr\}_{i,j=1,\ldots,k},
\end{eqnarray*}
and let the matrix $B \in \Rset^{k \times k}$ be defined by
\begin{eqnarray*}
&\displaystyle (F_a W_-)^{-1} F_a W_+ \e^{-(b-a) J_+} (F_b W_+)^{-1} F_b W_- &\\
 & \displaystyle   \bigl( I_{k \times k} - \e^{(b-a) J_-} (F_a W_-)^{-1} F_a W_+
\e^{-(b-a) J_+} (F_b W_+)^{-1} F_b W_-  \bigr)^{-1}.&
\end{eqnarray*}
Denoting by $\tau > 0$ the leading coefficient of $\Lop$, then
the integral $\int_a^b \Green_*(t,t)\, \df t$ equals the sum of the
following 8 terms:
\begin{eqnarray*}
\mathit{I} &=& N \tau^{-1} v_1 v_-, \\
\mathit{II} &=& \tau^{-1} v_1  \bigl( (F_a W_-)^{-1} F_a W_+ \odot A_{-+}
 \bigr) v_+, \\
\mathit{III} &=& -\tau^{-1} v_1  \bigl( (F_b W_+)^{-1} F_b W_- \odot A_{+-}
 \bigr) v_-, \\
\mathit{IV} &=& -\tau^{-1} v_1  \bigl( (F_b W_+)^{-1} F_b W_- \mathrm{e}^{(b-a)
  J_-} (F_a W_-)^{-1} F_a W_+ \odot A_{++}  \bigr) v_+, \\
\mathit{V} &=& \tau^{-1} v_1  ( B \odot A_{--}  ) v_-, \\
\mathit{VI} &=& \tau^{-1} v_1  \bigl( B \e^{(b-a)J_-} (F_a W_-)^{-1} F_a
W_+ \odot A_{-+}  \bigr) v_+, \\
\mathit{VII} &=& -\tau^{-1} v_1  \bigl( (F_b W_+)^{-1} F_b W_- \e^{(b-a)
  J_-} B \odot A_{+-}  \bigr) v_-, \\
\mathit{VIII} &=& -\tau^{-1} v_1  \bigl( (F_b W_+)^{-1} F_b W_- \mathrm{e}^{(b-a)
  J_-} B \e^{(b-a) J_-} (F_a W_-)^{-1} F_a W_+ \odot A_{++}  \bigr) v_+.
\end{eqnarray*}
\end{theo}

\begin{pf}
Since the characteristic polynomial has distinct roots
$\eta_1,\ldots,\eta_{2k}$, the Jordan canonical form of the
companion matrix is diagonal, and equation \eqref{eqGreen} implies that
$\partial_t^\mu (\Iop + \Delta\Lop)^{-1} \Eop_z(t)$ equals
\begin{eqnarray*}
&&\alpha_{2k}^{-1} \bar{v}_1 J^\mu \exp(t J) H^{-1} \bar{F}_a W \exp(a J)
 \biggl\{ \int_a^t \e^{-s \eta_i} \Eop_z(s)\, \df s
\cdot  ( W^{-1} v_2  )_i  \biggr\}_{i=1,\ldots,2k} \\[-2pt]
&& \quad {}- \alpha_{2k}^{-1} \bar{v}_1 J^\mu \exp(t J) H^{-1} \bar{F}_b W \exp(b J)
\biggl \{ \int_t^b \e^{-s \eta_i} \Eop_z(s)\, \df s
\cdot  ( W^{-1} v_2  )_i  \biggr\}_{i=1,\ldots,2k}.
\end{eqnarray*}
Since the function $\Eop_z$ is piecewise linear, the above
integrals can be explicitly evaluated over the intervals
$[t_j,t_{j+1}]$. For $j=0,N$, we have
\begin{eqnarray*}
  \int_a^{t_1} \e^{-s \eta_i} \Eop_z(s)\, \df s &=& \e^{-a
  \eta_i} \frac{1-\exp(-\Delta \eta_i / 2)}{\eta_i} z_1, \\[-2pt]
  \int_{t_N}^b \e^{-s \eta_i} \Eop_z(s)\, \df s &=& \e^{-t_N
  \eta_i} \frac{1-\exp(-\Delta \eta_i / 2)}{\eta_i} z_N,
\end{eqnarray*}
and for $j=1,\ldots,N-1$, we have
\begin{eqnarray*}
  &&\int_{t_j}^{t_{j+1}} \e^{-s \eta_i} \Eop_z(s)\, \df s\\[-2pt]
&& \quad = \int_0^\Delta \e^{-t_j \eta_i - \eta_i s}
   \bigl(  (1-s\Delta^{-1} ) z_j + s \Delta^{-1} z_{j+1}  \bigr)\,
  \df s \\[-2pt]
&& \quad = \e^{-t_j \eta_i} \int_0^\Delta \e^{-\eta_i s}
(1-\Delta^{-1} s)\, \df s  z_j
+ \e^{-t_j \eta_i} \int_0^\Delta \e^{-\eta_i s}  \Delta^{-1}
s\, \df s  z_{j+1} \\[-2pt]
&& \quad = \e^{-t_j \eta_i} \frac{\exp(-\Delta \eta_i) - 1 + \Delta
  \eta_i}{\Delta (\eta_i)^2} z_j
+ \e^{-t_j \eta_i} \frac{1 - (1 + \Delta \eta_i) \exp(-\Delta
  \eta_i)}{\Delta (\eta_i)^2} z_{j+1}.
\end{eqnarray*}
Arranging the eigenvalues as
$\eta_1^-,\ldots,\eta_k^-,\eta_1^+,\ldots,\eta_k^+$ and
inserting the definition of $\xi_-$, $\xi_+$, $\xi_-^0$, $\xi_+^0$,
$\xi_-^1$, $\xi_+^1$, we have that $\partial_t^\mu (\Iop +
\Delta\Lop)^{-1} \Eop_z(t_n)$ equals
\begin{eqnarray*}
 && 1_{n > 1} \sum_{j=1}^{n-1} \alpha_{2k}^{-1} \bar{v}_1 J^\mu \e^{t_n
  J} H^{-1} \bar{F}_a W \e^{(a-t_j) J}\pmatrix{ v_- \odot
  \e^{-\Delta J_-} \xi_-^0 \cr v_+ \odot \xi_+^0 } z_j \\[-2pt]
&& \quad {}+ \sum_{j=1}^n \alpha_{2k}^{-1} \bar{v}_1 J^\mu \e^{t_n J} H^{-1}
  \bar{F}_a W \e^{(a-t_{j-1}) J} \pmatrix{ v_- \odot (1_{j=1}
  \e^{-\Delta J_-/2} \xi_- + 1_{j > 1} \e^{-\Delta J_-} \xi_-^1) \cr
  v_+ \odot (1_{j=1} \xi_+ + 1_{j > 1} \xi_+^1) } z_j \\[-2pt]
&&  \quad   {}- \sum_{j=n}^N \alpha_{2k}^{-1} \bar{v}_1 J^\mu \e^{t_n J} H^{-1}
  \bar{F}_b W \e^{(b-t_j) J} \pmatrix{ v_- \odot (1_{j < N}
  \e^{-\Delta J_-} \xi_-^0 + 1_{j=N} \e^{-\Delta J_-/2} \xi_-) \cr v_+
  \odot (1_{j < N} \xi_+^0 + 1_{j=N} \xi_+)} z_j \\[-2pt]
&& \quad    {}- 1_{n < N} \sum_{j=n+1}^N \alpha_{2k}^{-1} \bar{v}_1 J^\mu \e^{t_n J}
  H^{-1} \bar{F}_b W  \e^{(b-t_{j-1}) J} \pmatrix{ v_- \odot
  \e^{-\Delta J_-} \xi_-^1 \cr v_+ \odot \xi_+^1 } z_j.
\end{eqnarray*}
The exponential factors on the terms $\xi_-$, $\xi_-^0$, $\xi_-^1$ may
be assimilated in the exponential factors before the large parenthesis
using $t_{j+1}-t_j = \Delta$ for $j=1,\ldots,N$ and
$t_2-t_1=t_{N+1}-t_N=\Delta/2$. Thereafter the terms in these sums are
of the same type as in equation \eqref{eqGreennum5} with $v_-$, $v_+$
replaced by  $v_- \odot \xi_-^0$, $v_+ \odot \xi_+^0$ etc., and the
formula for $\partial_t^\mu (\Iop + \Delta\Lop)^{-1} \Eop_z(t_n)$
follows by invoking the same reformulations as used in the proof of
Proposition~\ref{propnumericalstable}.

Finally, we consider the Green's function $\Green_*$ for $\Lop_* = v
\Iop + \Delta \Lop$. The differential operator $\Lop_*$ has leading
coefficient $\alpha_{2k} = \Delta \tau$, and inserting $s=t$ in the
first part of  equation \eqref{eqGreennum1}, we find that $\Green_*(t,t)$
equals
\begin{eqnarray*}
 & \displaystyle \Delta^{-1} \tau^{-1}  \bigl( v_1 \e^{-(b-t) J_-} - v_1
  \e^{-(b-t) J_+} (F_b W_+)^{-1} F_b W_-  \bigr) &\\[-2pt]
 &\displaystyle \bigl( I_{k \times k} - \e^{(b-a) J_-} (F_a W_-)^{-1} F_a W_+
\e^{-(b-a) J_+} (F_b W_+)^{-1} F_b W_-  \bigr)^{-1}& \\[-2pt]
&\displaystyle \bigl ( \e^{(b-t) J_-} v_- + \e^{(b-a) J_-} (F_a W_-)^{-1} F_a W_+
  \e^{-(t-a) J_+} v_+  \bigr).&
\end{eqnarray*}
To remove the possibly exploding exponential factor $\e^{-(b-t) J_-}$
in the first factor, we invoke the matrix formula
$(I-X)^{-1} = I + X (I-X)^{-1}$ on the second factor and rearranging
the exponential factors. Doing this $\Green_*(t,t)$ is rewritten
as the numerically stable expression
\begin{eqnarray*}
 &\displaystyle  \Delta^{-1} \tau^{-1}  \bigl( v_1 - v_1 \e^{-(b-t) J_+} (F_b
  W_+)^{-1} F_b W_- \e^{(b-t) J_-}  \bigr)& \\[-2pt]
 &\displaystyle \bigl( v_- + \e^{(t-a) J_-} (F_a W_-)^{-1} F_a W_+ \e^{-(t-a) J_+} v_+
   \bigr)& \\[-2pt]
&\displaystyle {}+ \Delta^{-1} \tau^{-1}  \bigl( v_1 \e^{(t-a) J_-} - v_1
  \e^{-(b-t) J_+} (F_b W_+)^{-1} F_b W_- \e^{(b-a) J_-}  \bigr)& \\[-2pt]
&\displaystyle (F_a W_-)^{-1} F_a W_+ \e^{-(b-a) J_+} (F_b W_+)^{-1} F_b W_-& \\[-2pt]
&\displaystyle\bigl ( I_{k \times k} - \e^{(b-a) J_-} (F_a W_-)^{-1} F_a W_+
 \e^{-(b-a) J_+} (F_b W_+)^{-1} F_b W_-  \bigr)^{-1}& \\[-2pt]
& \displaystyle \bigl( \e^{(b-t) J_-} v_- + \e^{(b-a) J_-} (F_a W_-)^{-1} F_a W_+
  \e^{-(t-a) J_+} v_+  \bigr).&
\end{eqnarray*}
This expression is expanded into the sum of 8 terms, which all may be
explicitly integrated over the interval $[a,b]$. For instance is the
integral over the second term given by
\[
  \int_a^b \Delta^{-1} \tau^{-1} v_1 \e^{(t-a) J_-} (F_a
  W_-)^{-1} F_a W_+ \e^{-(t-a) J_+} v_+\, \df t,
\]
which equals $\tau^{-1} v_1  ( (F_a W_-)^{-1} F_a W_+ \odot
A_{-+}  ) v_+$.
\end{pf}

\begin{remark*}
The predictors $\mv[x_m \vert y]$ may be seen as the predictors
$\mv[x^{\mathrm{fct}}_m \vert y]$ for the functional parameters
$x^{\mathrm{fct}}_m$ evaluated at the sample points $t_n$. The formulae
stated in Theorem~\ref{thmfda} may be extended to functional
representations for $\mv[x^{\mathrm{fct}}_m \vert y]$. Doing this the
predictions between sample points will be given as linear
combinations of exponential functions.
\end{remark*}

\begin{remark*}
If the kernel $\Green(t,s)$ is constant, say $\Green(t,s) = \lambda$,
then the operator approximation
\[
\int_0^1 \sum_{j=1}^N  ( v \Iop + \Mop_{\mathcal{T}}^{-1} \Lop
   )^{-1} \Eop_{e_j}(t_j)\, \df v =
  \int_0^1 \int_a^b \frac{\lambda \fracc{N}{b-a}}{1+N v \lambda}\, \df
  t\, \df v = \log(1 + N\lambda)
\]
gives the exact log determinant of $ \{ 1_{n=m} +
\Green(t_n,t_m)  \}_{n,m} = \Iop_N + \{ \lambda \}_{n,m}$. The
particular construction of the embedding operator $\Eop_z$ was chosen
to achieve this property.
\end{remark*}

A fundamental difference between our operator methods and the
smoothing spline technology lies in our dependence on boundary
conditions. Whether boundary conditions are desirable in statistical
modeling depends on the data situation at hand. If we have additional
knowledge implying particular boundary conditions, then this may be
used in the statistical model. However, in many data situations such
additional knowledge is not available, and the requirement to specify
boundary conditions may be disturbing. Here our advice is to use
Neumann-type conditions. Although the covariance function
$\Green(t,s)$ is not defined for Neumann conditions as noted in the
following example, this is possible due to the regularization induced
by the measurement noise; that is, $\Iop + \Mop_{\mathcal{T}}^{-1} \Lop$
is non-singular by construction.

\begin{example*} \label{examBrownian}
For $\Kop = \lambda \partial_t$ we have $\Lop = \Kop^\dag \Kop =
-\lambda^2 \partial_t^2$. Consider the following two sets of boundary
conditions:
\begin{eqnarray*}
(B1)\dvt \theta(a)  = \theta^{(1)}(b) = 0,  \qquad
(B2)\dvt \theta(a)  = \theta(b) = 0.
\end{eqnarray*}
We have $\Lop^{-1} \theta(t) = \int_a^b \Green(t,s) \theta(s)\, \df s$
with
\[
  \Green(t,s) = \cases{\displaystyle
\lambda^{-2}  \bigl( (t \wedge s)-a  \bigr), & \quad
\mbox{for boundary conditions (B1),} \cr\displaystyle
\lambda^{-2} \frac{((t \wedge s)-a)(b-(t \vee s))}{b-a}, & \quad
\mbox{for boundary conditions (B2).}
}
\]
Thus, the Laplace operator with boundary conditions (B1) leads
to the Brownian motion, and the Laplace operator with boundary
conditions (B2) leads to the Brownian bridge. The Laplace
operator with Neumann boundary conditions $\theta^{(1)}(a) =
\theta^{(1)}(b) = 0$ is not positive definite. Even so, this operator
can be used in a statistical model, where it implies an improper prior
for the serially correlated effects in terms of a Brownian motion with
a free level.

To compute the approximative log likelihood we find the Green's
function $\Green_v$ for $v \Iop + \frac{b-a}{N} \Lop$. In case of the
Brownian, motion equation \eqref{eqGreennum1} gives
{\fontsize{9.7pt}{12pt}\selectfont{
\begin{eqnarray*}
\Green_v(t,s)
&=& \frac{1}{\lambda} \sqrt{\frac{b-a}{N v}} \\
&&{}\times \frac{
\sinh ( (\fracf{(t \wedge s)-a}{\lambda \sqrt{b-a}} )\sqrt{N v}  )
\cosh ( (\fracf{b-(t \vee s)}{\lambda \sqrt{b-a}}) \sqrt{N v}  )}{
\cosh( \lambda^{-1} \sqrt{b-a} \sqrt{N v})},\\
\int_a^b \Green_v(t,t)\, \df t &=&
\frac{\sqrt{b-a}}{2\lambda} \sqrt{\frac{N}{v}} \cdot \frac{
\sinh( \lambda^{-1} \sqrt{b-a} \sqrt{N v})}{
\cosh( \lambda^{-1} \sqrt{b-a} \sqrt{N v})}.
\end{eqnarray*}}}%
In case of the Brownian bridge, equation \eqref{eqGreennum1}
gives
{\fontsize{9.7pt}{12pt}\selectfont{
\begin{eqnarray*}
  \Green_v(t,s)
&=& \frac{1}{\lambda} \sqrt{\frac{b-a}{N v}} \\
&&{}\times  \frac{
\sinh ( (\fracf{(t \wedge s)-a}{\lambda \sqrt{b-a}}) \sqrt{N v}  )
\sinh ( (\fracf{b-(t \vee s)}{\lambda \sqrt{b-a}}) \sqrt{N v}  )}{
\sinh( \lambda^{-1} \sqrt{b-a} \sqrt{N v})},\\
\int_a^b \Green_v(t,t)\, \df t &=&
\frac{\sqrt{b-a}}{2\lambda} \sqrt{\frac{N}{v}} \cdot \frac{
\cosh( \lambda^{-1} \sqrt{b-a} \sqrt{N v})}{
\sinh( \lambda^{-1} \sqrt{b-a} \sqrt{N v})} - \frac{1}{2v}.
\end{eqnarray*}}}%
In both cases the double integrals $\int_0^1 \int_a^b \Green_v(t,t)\,
\df t\, \df v$ can be computed giving explicit formulae for the
operator approximation of the matrix determinants. In case of an
equidistantly sampled Brownian motion, we have
\[
\log \det  \{ 1_{n=m} + \Green(t_n,t_m)  \}_{n,m=1,\ldots,N}
\approx \log \bigl( \cosh\bigl( \lambda^{-1} \sqrt{b-a} \sqrt{N} \bigr)  \bigr),
\]
and in case of an equidistantly sampled Brownian bridge, we have
\[
\log \det  \{ 1_{n=m} + \Green(t_n,t_m)  \}_{n,m=1,\ldots,N}
\approx \log\biggl ( \frac{ \sinh( \lambda^{-1} \sqrt{b-a} \sqrt{N} )}{
\lambda^{-1} \sqrt{b-a} \sqrt{N}}  \biggr).
\]
\end{example*}

\section{Approximative inference} \label{secstatistics}

In this section we combine the matrix formulae listed in
Section~\ref{secinference} with the operator approximation developed
in Section~\ref{secmath}. The obstacle in the matrix computations is the
inversion of the matrix $A_0 = \Iop_N + R_0 \in \Rset^{N \times N}$. Here
$R_0 = \{ \Green(t_n,t_m) \}_{n,m=1,\ldots,N}$ is defined via a
discretization $\mathcal{T} = \{t_1,\ldots,t_N\}$ and the Green's
function $\Green$ for a differential operator $\Lop = \sum_{l=1}^L
\Kop_l^\dag \Kop_l$.

The maximum likelihood estimator and the BLUPs given in
equations \eqref{eqMLE} and \eqref{eqBLUP} are approximated using the
block structure $A = A_0 \otimes \Iop_M$, the identity $A_0^{-1} z = z
- R_0 A_0^{-1} z$ for $z \in \Rset^N$ and the approximation
\[
R_0 A_0^{-1} z = A_0^{-1} R_0 z =  ( \Iop + R_0^{-1}  )^{-1} z \approx
   \{  ( \Iop + \Mop_{\mathcal{T}}^{-1} \Lop  )^{-1}
  \Eop_z(t_n)  \}_{n=1,\ldots,N}.
\]
Note that this approximation is applied both on the individual sample
vectors $y_m \in \Rset^N$ and on the sections of the columns of the
design matrices $\Gamma$ and $Z$. The approximation of the logarithmic
determinant equation \eqref{eqlogdet} in the restricted likelihood
equation \eqref{eqnlogRelikelihood} has already been stated in
equation \eqref{eqapprox2}, and the quadratic form of the serially
correlated effects is approximated by
\[
  \mv[x \vert y]^\top R^{-1} \mv[x \vert y] \approx
  \sum_{m=1}^M \sum_{n=1}^N \mv[x_m(t_n) \vert y]^\top
  \Mop_{\mathcal{T}}^{-1} \Lop \mv[x_m(t_n) \vert y].
\]
Furthermore, for an equidistant discretization with mesh length $\Delta$, we
have
\[
\mv[x \vert y]^\top R^{-1} \mv[x \vert y] \approx \Delta
  \sum_{l=1}^L \sum_{m=1}^M \sum_{n=1}^N  ( \Kop_l \mv[x_m(t_n) \vert y]
   )^\top  ( \Kop_l \mv[x_m(t_n) \vert y]  ).
\]

If the discretization $\mathcal{T}$ is equidistant, then semi-explicit
and numerically stable formulae for the above approximations are given
in Section~\ref{secexplicit}. For general discretizations the
operator approximations may be found as numerical solutions to
ordinary differential equations; for example, the function $f = (\Iop +
\Mop_{\mathcal{T}}^{-1} \Lop)^{-1} \Eop_z \in \mathcal{H}$ obeys to the
differential equation $f + \Mop_{\mathcal{T}}^{-1} \Lop f = \Eop_z$.

\section*{Acknowledgements}
I am grateful to Ib Skovgaard and Anders Tolver for encouragement
during the development of this work.


\printhistory

\end{document}